\documentclass[12pt,letterpaper]{article}
\usepackage{t1enc}
\usepackage[latin1]{inputenc}
\usepackage[english]{babel}
\usepackage{graphics}
\begin{document}

\title{Period function and characterizations of Isochronous potentials}
\author{A. Raouf Chouikha
\footnote
{Universite Paris 13 LAGA, Villetaneuse 93430.   
}}
\date{}
\maketitle

\begin{abstract} 
We are interested at first in the study of the monotonicity for the period function of the conservative equation \ $(1)\quad \ddot x + g(x) = 0.$\quad Some refinements of  known criteria are brought. Moreover, we give necessary and sufficient conditions so that the analytic potential of equation  $(1)$ is isochronous. These conditions which are different from those introduced firstly by Koukles and Piskounov and thereafter by Urabe appear sometime to be easier to use. We then apply these results to produce families of isochronous potentials depending on many parameters, some of them are news. Moreover, analytic isochronicity requirements of parametrized potentials will also be considered\\

\footnote
{2000 Mathematics Subject Classification  \ 34C05, 34C25, 34C35}. 
\end{abstract}

\bigskip
\section {Introduction and statement of results}

Consider the scalar equation with a center at the origin $0$
\begin{equation} \qquad \ddot x + g(x) = 0 \end{equation} or its planar equivalent system
\begin{equation} \qquad \dot x = y , \qquad \dot y = - g(x)\end{equation} 
where $\dot x = \frac{dx}{dt}, \ddot x = \frac{d^2x}{dt^2}$ and $g(x)$ is analytic on $R$.\\ Let $G(x)$ be the potential of equation (1) 
$$G(x) = \int_0^x g(\xi) d\xi.$$ We suppose in the sequel the following hypothesis holds

 $${(\cal H)} \quad \cases{  There \ exist\  a < 0 < b\ such \ that\ G(a) = G(b) = c,\ G(x) < c & \cr and \ x g(x) > 0 \ for \  all \ a < x < b \ and \ x \neq 0.\cr }$$ {\it Moreover, without lose of generality we will assume in the sequel that\\ $g(0) = 0$ and $g'(0) = 1.$}  \\
 
 Under these conditions system (2) admits a periodic orbit in the phase plane with energy $c$ and $a,b$ are the turning points of this orbit satisfying inequalities \ $\bar a < a < 0,\ 0 < b < \bar b,$ for some real $\bar a$ and $\bar b$.\\ This means the origin $0$ is a center of (2). This center is isochronous when the period of all orbits near $0 \in R^2$ are constant ($T = \frac {2\pi}{\sqrt{g'(0)}} = 2\pi$). The corresponding potential $G(x)$ is also called isochronous.  \\
 
Let $T(c)$ denotes the minimal period of this periodic orbit depending on the energy. It can be expressed under the following form
\begin{equation} T(c) = \sqrt{2}\int_a^b \frac{dx}{\sqrt{c-G(x)}}.\end{equation}  The period function is well defined for any $c$ such that $0 < c < \bar c$ and when $\bar c < +\infty$ then $\bar c$ is the energy of the homoclinic orbit with $\bar a < 0 < \bar b$ its turning points. It is well known that $T(c)$ is an analytic function of $c$.\\ 

Since the potential $G(x)$ has a local minimum at $0$, then we may consider an involution $A$ by 
$$ G(A(x)) = G(x) \ and \ A(x) x < 0 $$ for all $x \in [a,b]$. So, any closed orbit is $A$-invariant and $A$ exchanges the turning points: $b=A(a)$. \\ In fact, $A(x)$ is well defined in the interval $[a,b]$. To see that, set the function $$\rho(x)=\frac{x-A(x)}{2}.$$ This function is such that $\rho(A(x)) = - \rho(x)$ and $\rho'(x) = \frac{1-A'(x)}{2}$. Since $A'(x)<0$ we get $\rho'(x) >0$ and therefore $\rho$ is an analytic diffeomorphism on $[a,b]$. Then $A(x) = \rho^{-1}(-\rho(x))$ is well defined.  \\

Conversely, by using this involution we may calculate for a prescribed period function $T(c)$ the distance between the turning points. Indeed, Following Landau and Lifchitz we have [9, Chap.3, 12.1]): 
$$T(c) = \sqrt{2}\int_0^c [\frac{dA(a)}{dG}-\frac{d a}{dG}]\frac{dG}{\sqrt{c-G}}.$$
which implies 
\begin{equation} A(a) - a = b - a = \frac {1}{\pi \sqrt 2}\int_0^c \frac {T(\gamma) d\gamma}{\sqrt{c - \gamma}}.\end{equation}

 When the potential (or the center $0$) is isochronous or equivalently the period is constant for all orbits near zero  one has necessarily
$$b - a = 2 \sqrt{2c}.$$  

The behavior of the period function is important in applied mathematics. Many authors ([3],[4],[10],[11]) are interested in the monotonicity of this period function as well as in isochronicity cases of equation (1)  ([6],[7],[10]).\\ Isochronous centers plays a central role in several problems of dynamical systems. Also, it intervenes in the investigation of the  quantum spectrum. Many studies have been devoted to the problem of the relationship between  isochronicity equivalence and quantum isospectrality condition. See [1],[6] for more details. \\ 
In many papers one tried to characterize isochronous potentials. In the case of the rational potentials, it can be shown [2] that the only rational potential with a constant period which not reduces to a polynomial is the isotonic one 
$$G(x) = \frac{1}{8 \alpha ^2} [\alpha x +1 - \frac{1}{\alpha x + 1}]^2.$$

The paper is organized as follows :

- In the first part we study the monotonicity of the period function of a closed orbit of equation (1) depending on the energy and we produce better criteria than those known. In particular, we complete the work started by Chow and Wang [5].\\

- We propose to analyze the notion of isochronicity in using a new track. We will provide necessary and sufficient condition so that the center $0$ of system (2) be isochronous, or equivalently its analytic potential is isochronous. Many applications are brought. \\

- In the analytic case, let us write the potential under the following form $$ G(x)= \frac{1}{2}x^2+G_1(x)+G_2(x) $$ where $G_1(x)$ is odd and $G_2(x)$ is even. It is wellknown that the only symmetric isochronous potential is the harmonic one (i.e. $G(x)= \frac{1}{2}x^2$).  We prove here that when $G(x) $ is isochronous $G_1$ and $G_2$ are related  and consequently $G_2(x)\equiv 0$ implies necessarily $G_1(x)\equiv 0$.\\ Moreover, a simple proof of Chalykh-Veselov result [2] concerning rationals potentials will be provided. \\

- Finally, we derive a three-parameters family of isochronous potential more general than the one given by Dorignac [6]. This family includes harmonic and isotonic and Urabe potentials as well as the Bolotin-Mac Kay potential. Other families of isochronous potentials also will be considered.\\

\section {The period function}

Let us consider the period function of equation (1) depending on the energy \ $T\equiv T(c).$ \ When $g$ is an analytic function defined for\ $x \in [a,b],\quad T(c)$\ is also an analytic function defined for\ $c \in ]0,\bar c]$\ such that \ $lim_{c\rightarrow 0^+} T(c) = 2 \pi.$ \quad  
The following result states sufficient conditions for \ $T(c)$\ to be monotonic.

\bigskip

{\bf Theorem A} \quad {\it  Let $g(x)$ be an analytic function and $G(x) = \int_0^x g(\xi) d\xi$ be the potential of equation $(1)\quad \ddot x + g(x) = 0.$ Suppose hypothesis \ ($\cal{H}$)\ holds, let us define the n-polynomial with respect to $G$, $$f_n(G) = f(0)+f'(0)G+ \frac{1}{2}f''(0)G^2+....+\frac{1}{n!}f^{n}(0)G^n$$ where\quad  $f(0)= -g''(0)/3=-(1/3) \frac{d^2 g}{dx^2}(0),\quad f'(0)= -(7/9)g''^3(0)+ g^{(4)}(0)/5,\\ f''(0)= -\frac {1}{21}g^{(6)}(0) - \frac {310}{54} g''^5(0)+ 2g''^2(0) g^{(4)}(0),...$ \\ Suppose that for a fixed $n \in N$ and for $x \in [0,b]$ one has $${(\cal C}_n)\qquad \frac{d}{dx}[\frac{G}{g^2}(x)] > f_n(G) > \frac{d}{dx}[\frac{G}{g^2}(A(x))]$$ $$ (or\  \frac{d}{dx}[\frac{G}{g^2}(x)] < f_n(G) < \frac{d}{dx}[\frac{G}{g^2}(A(x))] ) $$ then the period function of (1)  $T'(c) > 0 \ (or < 0)$ for $0 < c < \bar c$.}\\ 
 
\bigskip

{\it Proof of Theorem A} \qquad  
Recall that 
$$  T(c) = \sqrt 2 \int_a^b \frac {dx}{\sqrt{c - G(x)}}.$$
In order to study the monotonicity of the period function $T(c)$ depending on the energy, it is sometime convenient to study its derivatives. 
We need the following \\

{\bf Lemma 2-1} \quad {\it The
 derivative of the period function (depending on the energy) $T'(c) = \frac{dT}{dc}$ may be written}
\begin{equation} T'(c) = (\frac {1}{c\sqrt 2})\int _a^b \frac {g^2(x) - 2 G(x) g'(x)}{g^2(x) \sqrt{c - G(x)} } dx\end{equation}

 This lemma has been initially proved by  Chow and Wang [5] but the proof we give below is different.\\
 
 {\it Proof} \qquad 
Consider the change $G(x) = s^2$, where $s=u(x)$ is a function of $x$ (this change has a sense because $G(x)$ is positive in a neighborhood of $0$ according to the hypothesis $xg(x) > 0$). It yields 
$$T(c) = 2\sqrt 2 \int_{-\sqrt c}^{\sqrt c} \frac {s ds}{g(u^{-1}(s))\sqrt{c - s^2}}.$$ 
By another change $s= {\sqrt c}\sin\theta$ the period function may be expressed 
$$T(c) = \sqrt 2 \int_{-\pi/2}^{\pi/2 } \frac { d\theta}{u'(u^{-1}(\sqrt c)\sin\theta))} = \sqrt 2 \int_{-\pi/2}^{\pi/2} (u^{-1})'(\sqrt c \sin\theta) d\theta .$$ Then we get by another way the derivative of the period function
$$T'(c) = \frac {2}{\sqrt c} \int_{-\pi/2}^{\pi/2} (u^{-1})''(\sqrt c\sin\theta)\sin\theta d\theta .$$ 
By splitting $T'(c) = \frac {2}{\sqrt c}\int_{-\pi/2}^0 \quad +\quad \frac {2}{\sqrt c}\int_0^{\pi/2}$ and using the formula\\ $(u^{-1})''(u(s)) = -\frac{u''(s)}{(u'(s))^3}$ and making another time the change of variable $s=u(x)= {\sqrt c}\sin\theta$ we then obtain 
$$ T'(c) = \frac {1}{c\sqrt 2}\int_a^0 \frac {g^2(x) - 2 G(x) g'(x)}{g^3(x)}\frac{ -g(x)}{\sqrt{c - G(x)} } dx + $$ $$\frac {1}{c\sqrt 2}\int_0^b \frac {g^2(x) - 2 G(x) g'(x)}{g^3(x)}\frac{ g(x)}{\sqrt{c - G(x)} } dx.$$ Since $x g(x) <0$ the derivative can be written  $$T'(c) = (\frac {1}{c\sqrt 2})\int _a^b \frac {g^2(x) - 2 G(x) g'(x)}{g^2(x) \sqrt{c - G(x)} } dx.$$
\bigskip

To prove Theorem A observe at first   
\begin{equation} T'(c) = (\frac {1}{c\sqrt 2})\int _a^b \frac {\frac{d}{dx}(G/g^2)}{ \sqrt{c - G(x)} } g(x)dx.\end{equation}
Write $$(c\sqrt 2)T'(c) = \int _a^0 \frac {g^2(x) - 2 G(x) g'(x)}{g^2(x) \sqrt{c - G(x)} } dx + \int _0^b \frac {g^2(x) - 2 G(x) g'(x)}{g^2(x) \sqrt{c - G(x)} } dx.$$  By using the involution $A$ and since $x(g(x) > 0$ for $x\neq 0$ we get 
$$T'(c) =   (\frac {1}{c\sqrt 2})\int _b^0 \frac {\frac{d}{dx}(G/g^2)(A(x))}{ \sqrt{c - G(x)} } g(A(x))A'(x)dx + (\frac {1}{c\sqrt 2})\int _0^b \frac {\frac{d}{dx}(G/g^2)}{ \sqrt{c - G(x)} } g(x)dx$$
$$ =   (\frac {1}{c\sqrt 2})\int _b^0 \frac {\frac{d}{dx}(G/g^2)(A(x))}{ \sqrt{c - G(x)} } g(x)dx + (\frac {1}{c\sqrt 2})\int _0^b \frac {\frac{d}{dx}(G/g^2)}{ \sqrt{c - G(x)} } g(x)dx$$
$$= (\frac {1}{c\sqrt 2}) \int_0^b \frac{\frac{d}{dx}(\frac{G}{g^2}(x))-\frac{d}{dx}(\frac{G}{g^2}(A(x))}{\sqrt{c-G(x)}}g(x) dx.$$

 Then 
$$T'(c) > (or\ <)\ (\frac {1}{c\sqrt 2})\int _a^b \frac {f(0)+f'(0)G+ \frac{1}{2}f''(0)G^2+....+\frac{1}{n!}f^{n}(0)G^n}{ \sqrt{c - G(x)} } g(x)dx .  $$
Let us integrate by parts the following for $p \geq 1$
$$\int \frac {G^p}{ \sqrt{c - G(x)} } g(x)dx =\int \frac {G^p}{ \sqrt{c - G(x)} } dG = -G^p \sqrt{c-G} + 2p \int G^{p-1} \sqrt{c-G} dG.$$
Therefore, it is easy to see that the following integral should be zero since $G(a) = G(b) =c$ $$\int _a^b \frac {G^p}{ \sqrt{c - G(x)} } g(x)dx = 0$$ for any integer $0 \leq p \leq n.$\\
Thus, the right side of (6) is non zero and condition $$\frac{d}{dx}[G(x)/g^2(x)] - f_n(G) \neq 0 $$ implies $T'(c)\neq 0$. This means \ $T=T(c)$ \ is monotonic.\\

\bigskip 

As applications of Theorem A, we derive at first monotonicity conditions for the period function $T(c)$ depending on the energy. This problem holds  importance in dynamical systems (see [4], [10] and [11] for example). In this part, we find again some known criteria for the monotonicity of the period function and propose some others which are better and seem to be new in the literature.\\  

For a complete study and a comparison between these sufficients conditions we refeer to [3] and [4] and references therein.\\ 
Notice that the monotonicity criteria for the period function produced by Cor. 2.5 of [5] is more general than those given by F. Rothe [10] and R. Schaaf [11].\\
We propose the following which slightly improves Cor. 2.5 of [5].

\bigskip
{\bf Corollary 2-2}\quad {\it Suppose hypothesis ($\cal{H}$) holds and let $g(x)$ be an analytic function for $ \bar a < x < \bar b $ and \\ $G(x) = \int_0^x g(\xi) d\xi$ be the potential of (1) and $g''(0) \neq 0$.\\
1 - Suppose condition ${(\cal C}_0)$ holds, this means 
$$g^2(x) + \frac {g''(0)}{3} g^3(x) - 2 G(x) g'(x) > 0 (or < 0)$$  then the period function of (1) is such that $T'(c) > 0 (or < 0)$ for $0 < c < \bar c$.\\
2 - Suppose \begin{equation} \qquad \frac{d}{dx}[G(x)/g^2(x)] = \alpha \end{equation} 
where $\alpha = - g''(0)/3 \neq 0$  then $0$ is an isochronous center of $(1)$.\\ 
Moreover, the so-called Urabe potential 
$$G(x) = \frac {4}{\alpha^2} - \frac {2}{\alpha}(x + 2\frac {\sqrt{1-\alpha x}}{\alpha})$$
where  $\mid x\mid < 1/\alpha$ is the unique analytic solution of (7) verifying \\ $ g(0) = G(0) = 0$ and $G''(0)=1.$}\\

\bigskip
{\it Proof of Corollary 2-1}\quad 
Indeed,  (7)  may be written $$ g^2(x) + \frac {g''(0)}{3} g^3(x) - 2 G(x) g'(x) = g^3(x)\frac{d}{dx}([G(x)/g^2(x)]+ \frac {g''(0)}{3} g^3(x) .$$ Then (7) must be  
$$\frac{d}{dx}[G(x)/g^2(x)] = -g''(0)/3 = \alpha $$
Implying $$G(x)= \frac{1}{2} g^2(x)+ \alpha x g^2(x)$$  since \ $ g'(0)=1$. Thus, starting from (6) the derivative of the period function $T(c)$ may be expressed under the following forms 
$$T'(c) = (\frac {1}{c\sqrt 2})\int _a^b \frac {g^2(x) - 2 G(x) g'(x)}{g^2(x) \sqrt{c - G(x)} } dx$$ 
$$ = (\frac {1}{c\sqrt 2})\int _a^b \frac {\frac{d}{dx}(G/g^2)}{ \sqrt{c - G(x)} } g(x)dx$$
$$ = (\frac {1}{c\sqrt 2})\int _a^b \frac {\frac{d}{dx}({\sqrt G}/g)^2}{ \sqrt{c - G(x)} } g(x)dx.$$ Integrating by parts, one gets
$$T'(c) = (\frac {1}{c\sqrt 2})\int _a^b \frac {\alpha }{ \sqrt{c - G(x)} } g(x)dx \equiv 0.$$

  On the other hand, since
$$\frac{g^2-2 G g'}{g^3}
 = \frac{d}{dx}(G/g^2)= \frac{d}{dx}({\sqrt G}/g)^2\qquad  $$ it implies that $G=G(x)$ may be inverted and that $x$ can be expressed 
 $$x = \sqrt{2G} + \alpha G.\qquad (7')$$
 
 {\bf Remarks 2-3} \ Inversing now equation (7') and according to initial conditions \ $ g(0) = G(0) = 0$ and $G''(0)=1$ \ it is easy to see that the potential
 $$G(x)=\frac {4}{\alpha^2} - \frac {2}{\alpha}(x + 2\frac {\sqrt{1-\alpha x}}{\alpha})$$
 is the unique solution of $(7')$ or $(7)$. $G(x)$ is called {\it Urabe potential}.\\ Notice that its derivative $g(x)$ must also be a solution of 
 $$5 (g''(x))^2 - 3 g'(x) g'''(x) = 0$$ and consequently for any Urabe potential we get necessarily 
\begin{equation} g^{(4)}(0) = \frac{35}{9}g''^3(0).\end{equation}

 When in addition $g''(0)=0$ (i.e. $\alpha = 0$) equation (7) has as unique solution the harmonic potential $G(x) = (1/2) x^2$ .\\  
In fact, as we will see below there exist isochronous potentials that do not verify the restrictive condition (8). Consider the following expression 
\begin{equation} {\cal F}(x) = g^2(x) + \frac {g''(0)}{3} g^3(x) - 2 G(x) g'(x) + (\frac {7g''^3(0)}{9}- \frac {g^{(4)}(0)}{5}) g^3(x) G(x).\end{equation} 
The following brings a new monotonicity condition for the period function \\ 

\bigskip

{\bf Corollary 2-4}\quad {\it Let $g(x)$ be an analytic function  and \\ $G(x) = \int_0^x g(\xi) d\xi$ be the potential of equation (1). Suppose ${(\cal C}_1)$ holds, this means  $$g^2(x) + \frac {g''(0)}{3 } g^3(x) - 2 G(x) g'(x) + (\frac {7g''^3(0)}{9}- \frac {g^{(4)}(0)}{5}) g^3(x) G(x) > 0 (or < 0)$$ for $ \bar a < x < \bar b $ then $T'(c) > 0 (or < 0)$ for $0 < c < \bar c$. So a sufficient condition for (1) to have an isochronous center is $${\cal F}(x) = g^2(x) + \frac {g''(0)}{3} g^3(x) - 2 G(x) g'(x) + (\frac {7g''^3(0)}{9}- \frac {g^{(4)}(0)}{5}) g^3(x) G(x) = 0$$ or equivalently  $$\frac{d}{dx}[G(x)/g^2(x)] = \alpha + \beta G $$ where $\alpha = -\frac {g''(0)}{3}$ and $\beta =- (\frac {7g''^3(0)}{9}- \frac {g^{(4)}(0)}{5}) \neq 0$. }\\
\bigskip

 {\it Proof of Corollary 2-4}\quad 
Indeed, recall that $$\frac{g^2-2 G g'}{g^3}
 = \frac{d}{dx}(G/g^2)= \frac{d}{dx}({\sqrt G}/g)^2 $$
 Then ${\cal F}(x)\equiv 0$ means that 
$$\frac{g^2-2 G g'}{g^3} = \alpha + \beta G $$ 
where $\alpha = -g''(0)/3$ and $\beta= -(7/9)a^3 + \gamma/5 = -(7/9)g''^3(0)+ g^{(4)}(0)/5$.\\
Consider again the derivative of the period function 
$$T'(c) =   (\frac {1}{c\sqrt 2})\int _a^b \frac {\frac{d}{dx}(G/g^2)}{ \sqrt{c - G(x)} } g(x)dx =(\frac {1}{c\sqrt 2})\int _a^b \frac {\frac{d}{dx}({\sqrt G}/g)^2}{ \sqrt{c - G(x)} } g(x)dx.$$ So one obtains
$$T'(c) = (\frac {1}{c\sqrt 2})\int _a^b \frac {\alpha + \beta G}{ \sqrt{c - G(x)} } g(x)dx$$
$$=(\frac {1}{c\sqrt 2})\int _a^b \frac {\alpha}{ \sqrt{c - G(x)} } g(x)dx + (\frac {1}{c\sqrt 2})\int _a^b \frac { \beta G}{ \sqrt{c - G(x)} } g(x)dx.$$
By Corollary 2-2, the first integral of the right side should be $0$. The second integral may be written as 
  $$(\frac {1}{c\sqrt 2})\int _a^b \frac { \beta G}{ \sqrt{c - G(x)} } g(x)dx = (\frac {1}{c\sqrt 2})\int _a^b \frac { \beta (G-c)}{ \sqrt{c - G(x)} } g(x)dx $$ $$= (\frac {1}{c\sqrt 2})\int _a^b  \beta  \sqrt{c - G(x)} g(x)dx = (\frac {1}{c\sqrt 2})\int _{G(a)}^{G(b)}  \beta  \sqrt{c - G} dG = 0$$
Thus, one gets
$$T'(c) \equiv 0.$$

Moreover, we also may deduce
$$2\frac{d}{dx}({\sqrt G}/g) = \beta g {\sqrt G} + \alpha g/{\sqrt G}.$$

On the other hands $G=G(x)$ may be inverted and $x$ can be expressed in terms of $G$.
By integration and according to hypotheses $G(0) =g(0) = 0$ and $g'(0)=1,$ it yields
$$2 ({\sqrt G}/g) = (2\beta /3) G{\sqrt G} + 2\alpha {\sqrt G}+ {\sqrt 2}, $$
or equivalently
$$1 = (\beta /3) g G + \alpha g + (\frac{g}{\sqrt {2G}}) $$
So by another integration one obtains 
$$x = (\beta /6) G^2 + \alpha G + \sqrt {2G} $$
$$x = \sqrt {2G} - \frac {g''(0)}{3}G + (\frac {7g''^3(0)}{54}- \frac {g^{(4)}(0)}{30})G^2.$$

{\bf Remarks 2-5}\quad Notice that clearly Corollary 2-4 is more general than Corollary 2-2. We then obtain a better criteria of the monotonicity for the period function. In the sense that ${(\cal C}_0)$ implies ${(\cal C}_1)$ which implies $T'(c) > 0 $ (or $< 0$).\\ By the same way for any fixed $n$ many other sufficient conditions ensuring the monotonicity of the period function of the form 
$$ g^2(x) -2G(x) g'(x) - \alpha g^3(x) - \beta g^3(x) G(x) - \gamma g^3(x) G^2(x)- .....> 0 \ (or\ <0)$$
may be deduced. In particular, condition $$ g^2(x) -2G(x) g'(x) - \alpha g^3(x) - \beta g^3(x) G(x) - \gamma g^3(x) G^2(x)- .....= 0$$
implies that the potential $G(x)$ is isochronous and $x$ can be expressed $$x =   \sqrt {2G}+ \alpha G + (\beta /6) G^2 + (\gamma/24) G^3+....$$
Thus, for a n-polynomial $f_n(G) = f(0)+f'(0)G+ \frac{1}{2}f''(0)G^2+....+\frac{1}{n!}f^{n}(0)G^n$  conditions of Theorem A
$${(\cal C}_n)\qquad \frac{d}{dx}[\frac{G}{g^2}(x)] > ( <) f_n(G) >( <) \frac{d}{dx}[\frac{G}{g^2}(A(x))]$$  ensuring $T'(c) > 0 \ (or\ < 0)$ are such that ${(\cal C}_{n+1})$ is finest than ${(\cal C}_n)$.

 \section{Isochronicity conditions for a center of equation (1)}
 
 \subsection{Isochronicity conditions for a center of equation (1)}

The problem to determine whether the center is isochronous has attracted many researchers for long time. This problem has been recently revived due to advancement of computer algebra. New powerfull algorithms have been discovered indeed. \\ For the sake of completeness let us recall below different criteria for the isochronicity of periodical solutions of equation (1).\\

Using formula (4), Landau and Lifschitz [9] deduced the following

\bigskip 

 {\bf Proposition 3-1} \ [9, Chap.3]  When $g(x)$ is continuous and hypothesis ${(\cal H)}$ holds, (1) has an isochronous center at the origin $0$ if and only if 
 $$x - A(x) = 2\sqrt{2 G(x)}$$
 for all $0<x<b$ where $A$ is the involution such that $G(A(x)=G(x)$ and $A(x) x < 0$ for $x \neq 0$. 

\bigskip

When $g(x)$ is a continuous function, $g(x)$ and $x$ having the same sign,  Koukles and Piskounov [7] produced necessary and sufficient conditions so that the center of the system (2) is isochronous. \\ 

\bigskip
 {\bf Proposition 3-2} \ [7, Th 5] \quad {\it A set of necessary and sufficient conditions for the period of every solution of (1) near $0$ to be equal to a constant $T_0$ is:\\
 1 - $g(x)$ is continuous and positive for small positive $x$.\\
 2 - $liminf_{x\rightarrow 0} \mid\frac{g(x)}{x}\mid \neq 0.$\\
 3 - $T_0 \geq limsup_{x\rightarrow 0} \frac{2\pi}{\sqrt{g(x)/x}}.$\\
 4 - $g(-x) = - \frac{d}{dx}[\frac{T_0}{\pi}\sqrt{2x}-(\int_0^xg(u) du)^{-1}]^{-1}$\ where index $-1$ denotes an inverse function.}\\

\bigskip
This result was improved by Koukles and Piskounouv themselves in the analytic case. They proved the following \\

\bigskip
 {\bf Proposition 3-3} \ [7, Th 6] \quad {\it Let $g(x)$ be a real analytic function.  Then the center $0$ of the equation $(1)\quad \ddot x+g(x) = 0$ is isochronous if and only if the inverse function $x = \Theta(z)$ \ 
 $$ [\int_0^x g(\xi) d\xi ]^{-1} = \Theta(z)$$ is of the form
 $$ \Theta(z) = \sqrt{z}+ P(z)$$ where $P$ is a real analytic function such that $P(0) =0$ . }
 
  \bigskip
  Later, Urabe proposed some refinements of Proposition 3-3 by considering the assumption of the differentiability of $g(x)$ at $0$ and proved the following result which is most used than Propositions 3-2 and 3-3
   
\bigskip
 
 {\bf Proposition 3-4} \ [13] \quad {\it Let \ $ g(x) $\ be a \ $C^1$ \  function defined in $V_0$ a neighborhood of $0$ verifying   $x g(x) > 0$ in $V_0 /\{0\}$. Then the system $(2)$ has an isochronous center at the origin $0$ if and only if   \  $g(x)$\ may be written 
 $$g(x) = \frac {X}{1 + h(X)}$$ 
 where \ $h(X)$\ is a \ $C^1$ \ odd function and \ $X = \sqrt {2 G(x)},\ \frac {X}{x} > 0$\ for \ $x\neq 0$.}
 \bigskip

We will propose an alternative approach in order to derive isochronous potentials. Some other criteria or equivalent characterizations will be presented. Their significance makes the study of the isochronicity much easier, since any isochronous potential $G$ appears to be solution of a differential equation. More precisely, we state the following

\bigskip 

{\bf Theorem B} \quad {\it  Suppose hypothesis ${(\cal H)}$ holds and let $g(x)$ be an analytic function and $G(x) =\int_0^x g(s) ds $.\ Then the  equation $$ \ddot x + g(x) = 0 \qquad (1)$$  has an isochronous center at $0$ if and only if  
\begin{equation}\frac{d}{dx}[G(x)/g^2(x)] = f(G) \end{equation} 
where $f$ is an analytic function defined in some neighborhood of $0$. }

\bigskip

We first deduce from Theorem B that the analytic functions $f$ and $g$ are naturally related and by (10) one has necessarily $$ f(0)= -g''(0)/3, f'(0)= -(7/9)g''^3(0)+ g^{(4)}(0)/5,$$  $$f''(0)= -\frac {1}{21}g^{(6)}(0) - \frac {310}{54} g''^5(0)+ 2g''^2(0) g^{(4)}(0),....$$\\

\bigskip

{\it Proof of Theorem B} \quad 
Recall that 
$$  T(c) = \sqrt 2 \int_a^b \frac {dx}{\sqrt{c - G(x)}}$$
 and by Lemma 2-1 its derivative $T'(c) = \frac{dT}{dc}$ may be written
\begin{equation} T'(c) = (\frac {1}{c\sqrt 2})\int _a^b \frac {g^2(x) - 2 G(x) g'(x)}{g^2(x) \sqrt{c - G(x)} } dx\end{equation}
\begin{equation} = (\frac {1}{c\sqrt 2})\int _a^b \frac {\frac{d}{dx}(G/g^2)}{ \sqrt{c - G(x)} } g(x)dx.\end{equation}


Suppose at first that $0\in R^2$ is an isochronous center of system (2). Then, by Proposition 3-1 
\begin{equation}
x - A(x) = 2\sqrt{2G(x)}\qquad for \quad all \qquad 0<x<b.
\end{equation} 
Deriving this expression 
$$1-A'(x) = \frac{g(x)\sqrt 2}{\sqrt G(x)}$$
it implies
$$\frac{G(x)}{g^2(x)} = \frac{2}{(1-A(x)^2}.$$
Therefore, by deriving with respect to $x$
\begin{equation}
\frac{d}{dx}[\frac{G(x)}{g^2(x)}] = \frac{d}{dx}[\frac{2}{(1-A(x)^2}] = \frac{4 A''(x)}{(1-A'(x))^3}
\end{equation}
On the other hand, we need the following

\bigskip
 {\bf Lemma 3-1}\quad {\it For any analytic involution $A(x)$ defined for all $x \in [a,b]$ the following expression holds }
$$\frac{ A''(A(x))}{(1-A'(A(x)))^3}=\frac{ A''(x)}{(1-A'(x))^3}$$ 
\bigskip 
 {\it Proof} \quad 
Indeed, to prove this lemma 
we derive $A(A(x)) = x$ \ it yields $$A'(A(x)) = \frac{1}{A'(x)}\qquad and \qquad A''((A(x)) = -\frac{A''(x)}{A'^3(x)}.$$ 
We now replace $x$ by $A(x)$ in the right side of (9), one obtains
$$\frac{ A''(A(x))}{(1-A'(A(x)))^3}= - \frac{A''(x)}{A'^3(x)}\frac{1}{(1-\frac{1}{A'(x)})^3}$$
So, since $A'(x) \neq 0$ one gets after simpification
$$\frac{4 A''(A(x))}{(1-A'(A(x)))^3}=\frac{4 A''(x)}{(1-A'(x))^3}$$

This lemma implies that expression (14) is $A$-invariant and therefore $\frac{d}{dx}[\frac{G(x)}{g^2(x)}] $ is an analytic function only dependent on  $G$ since $G(A(x))= G(x)$.\\

  We have now to prove the converse. We will use for that the following \\

  {\bf Lemma 3-2}\quad {\it The derivative of the period function $T(c)$ may be written as}
$$T'(c) = (\frac {1}{c\sqrt 2}) \int_0^{b}\frac{\frac{d}{dx}(\frac{G}{g^2}(x))-\frac{d}{dx}(\frac{G}{g^2}(A(x))}{\sqrt{c-G(x)}}g(x) dx.$$
\bigskip

{\it Proof}\quad Indeed, we have seen that $$T'(c) = (\frac {1}{c\sqrt 2})\int _a^b \frac {\frac{d}{dx}(G/g^2)(x)}{ \sqrt{c - G(x)} } g(x)dx$$ for $\bar a < a < 0$ and $0  < b < \bar b$. \\ By splitting the integral we get
$$T'(c) =   (\frac {1}{c\sqrt 2})\int _b^0 \frac {\frac{d}{dx}(G/g^2)(A(x))}{ \sqrt{c - G(x)} } g(A(x))A'(x)dx + (\frac {1}{c\sqrt 2})\int _0^b \frac {\frac{d}{dx}(G/g^2)}{ \sqrt{c - G(x)} } g(x)dx$$
$$ =   (\frac {1}{c\sqrt 2})\int _b^0 \frac {\frac{d}{dx}(G/g^2)(A(x))}{ \sqrt{c - G(x)} } g(x)dx + (\frac {1}{c\sqrt 2})\int _0^b \frac {\frac{d}{dx}(G/g^2)}{ \sqrt{c - G(x)} } g(x)dx$$
$$= (\frac {1}{c\sqrt 2}) \int_0^b \frac{\frac{d}{dx}(\frac{G}{g^2}(x))-\frac{d}{dx}(\frac{G}{g^2}(A(x))}{\sqrt{c-G(x)}}g(x) dx.$$

Thus, by this lemma condition $(10) \quad \frac{d}{dx}[G(x)/g^2(x)] = f(G) $ implies necessarily that $\frac{d}{dx}[G(x)/g^2(x)] = \frac{d}{dx}[G(x)/g(A(x))^2]) $ and then $T'(c) \equiv 0$.\\
Theorem B is then proved.\\

 Let us consider $F$ the primitive of $f $ such that $F(0)=0 $. Then \\ $$(10) \ \frac{d}{dx}[G(x)/g^2(x)] = f(G(x)) \Leftrightarrow g(x)\frac{d}{dx}[G(x)/g^2(x)] = g(x)f(G(x)).$$  Integrate by parts, it yields \ $ \frac{G(x)}{g(x)}- x = F(G(x)).$ \ Therefore (10) is equivalent to \ $2G(x) - xg(x) = g(x) F(G(x)) $.\ Thus, we get another criteria of isochronicity

\bigskip 

{\bf Corollary 3-3} \quad {\it Under hypotheses of Theorem A, equation (1) admits an isochronous center at $0$ if and only if 
$$2G(x) - xg(x) = g(x) F(G(x)) \qquad (10')$$
where $F$ is an analytic function defined in some neighborhood of $0$. }

\bigskip 
 
 \bigskip
 
Some other equivalent conditions may also be deduced

\bigskip
{\bf Corollary 3-4} \quad {\it Under hypothesis $({\cal H}),\ 0$ is an isochronous center of (1) if and only if $x=x(G)$ is an analytic solution of the linear ODE 
\begin{equation} 2G \frac{d^2x}{dG^2}+\frac{dx}{dG} = f(G),\end{equation}
where $f$ is an analytic function. \\ Moreover, this solution must satisfy the conditions:} $$x(0)= 0,\quad lim_{G\rightarrow 0}(\frac{x^2}{2G})=1.$$\\

\bigskip

{\it Proof of Corollary 3-4}\quad Let us consider again $F(t)$ the integral of $f(t)$. Then, it is easy to see that condition \ $\frac{d}{dx}(G/g^2) = f(G)$\ is equivalent to \begin{equation}2G \frac{dx}{dG}= x+ F(G).\end{equation}
We derive (16) with respect to the variable $G$, we then obtain (15).\\
So, for any analytic function $f$ the linear equation (15) admits a unique solution $x = x(G)$ according to initial conditions. More precisely, consider the change \ $x = \sqrt{2G} + y.$\ Then \ $y=y(G)$ is solution of $$ 2G \frac{d^2y}{dG^2}+\frac{dy}{dG} = f(G)$$ with initial points : $y(0) =0, y'(0)= 1.$  A resolution of the last equation yields 
$$y(G)= \sqrt{2G}\int_0^G \frac{F(\nu)}{(2\nu)^{3/2}}d \nu.$$ Thus, a solution of (15) may be written
$$x(G)= \sqrt{2G}\ (1+ \int_0^G \frac{F(\nu)}{(2\nu)^{3/2}}d \nu).$$ \\

Let us consider now the involution $A$ such that \ $G(A(x))= G(x)$.\  Another consequence of Theorem B is
\bigskip

{\bf Corollary 3-5} \quad {\it  Let $G(x) =\int_0^x g(s) ds $ be an analytic potential of the scalar equation $$ \ddot x + g(x)= 0. \qquad (1)$$ \\ Let $A(x)$ be  an analytic involution defined by: \\ $  
G(A(x)) = G(x)$ \ and \ $A(x)x < 0$.\\ Then $0$ is an isochronous center of (1) if and only if  $$\frac{d}{dx}[G(x)/g^2(x)] = \frac{4 A''(x)}{(1-A'(x))^3}. $$   Moreover, the last expression is $A$-invariant. }\\

\bigskip
{\it Proof of Corollary 3-5}\quad Let $A$ be an analytic involution then \\ $G(A(x)) = G(x)$\ implies 
$$\frac{dG}{dx}= g(x) = \frac{dA}{dx} g(A(x))$$ and $$\frac{G}{g^2}(A(x)) = (\frac{dA}{dx})^2 \frac{G}{g^2}(x).$$ Deriving the last expression we then  obtain $$A'(x) \frac{d(\frac{G}{g^2})}{dx}(A(x))=2 A'(x) A''(x) \frac{G}{g^2}(x)+A'^2(x) \frac{d(\frac{G}{g^2})}{dx}(x).$$ Since $A'(x) \neq 0$ then, after simplification we get the differential equation
\begin{equation} \frac{d(\frac{G}{g^2})}{dx}(A(x))=2  A''(x) \frac{G}{g^2}(x)+A'(x) \frac{d(\frac{G}{g^2})}{dx}(x). \end{equation}
By Theorem B,\\ $\frac{d}{dx}[G(x)/g^2(x)] = f(G)$ implies $\frac{d}{dx}[G(x)/g^2(x)] = \frac{d}{dx}[G(x)/g(A(x))^2]) .$\\ Therefore, the solution of  equation (17) is \begin{equation} \frac{G}{g^2}(x) = \frac{2}{(1-A'(x))^2}\end{equation} since  $A'(0) = -1.$ Thus, by deriving one gets \  $\frac{d}{dx}[G(x)/g^2(x)] = \frac{4 A''(x)}{(1-A'(x))^3}. $\\
To prove the converse we require again Lemma 3-2 

$$T'(c) = (\frac {1}{c\sqrt 2}) \int_0^{b}\frac{\frac{d}{dx}(\frac{G}{g^2}(x))-\frac{d}{dx}(\frac{G}{g^2}(A(x))}{\sqrt{c-G(x)}}g(x) dx.$$



\bigskip 

The analytic involution $A(x)$ can also be defined as a solution of a linear ODE. The following result is analogous to Corollary 3-4.

\bigskip 
{\bf Corollary 3-6} \quad {\it Under hypotheses of Theorem A, equation (1) admits an isochronous center at $0$ if and only if the involution $A=A(G)$ is a solution of 
\begin{equation} 2G \frac{dA}{dG}= A(G) + F(G)
\end{equation}
 $f$ is an analytic function and $F$ is its integral.\\ Moreover, this solution must satisfy the conditions:} $$A(0)= 0,\quad lim_{G\rightarrow 0}(\frac{A^2}{2G})=1$$

\bigskip 

{\it Proof of Corollary 3-6}\quad To see that, we start from  $A(A(x))= x$\ and \ $G(A(x) = G(x). $ \ So equation (11) \ $2G \frac{dx}{dG}= x + F(G)$\ is equivalent to 
$$2G(A(x)= \frac{dG(A(x))}{dx}[ A(x) + F(G(A(x))]= \frac{dG(x)}{dx}(x + F(G)) .$$ On the other hand, deriving  $G(A(x) = G(x) $ with respect to $x$ we get
\ $ (\frac{dA(x)}{dx})(\frac{dG(A(x))}{dx})= \frac{dG(x)}{dx}$.\ We then obtain $$A(x) + F(G(x)) = [x + F(G(x))]\frac{dA(x)}{dx}.$$ Hence we deduce the ODE (19).\\

{\bf Remark 3-7} \quad The method of Urabe requires the use an intermediary function $h$ that is not in general explicitly known.  Indeed, in order to determine the potential isochrones one must be able to show by means of the change of variables $X=\sqrt{2G}$ that the inverse function of $G(x) =y $ is of the form $x = X+ H(X) $. What is not easy to achieve in any case. On the other hand, our approach is more direct. Our criteria are simply relations between an isochronous potential and its derivative. To be clearer, condition (10) of Theorem B is a sort of non linear differential equation of order one with respect to the variable $x$. All solution $G=G(x)$ of $(10) \ \frac{d}{dx}[G(x)/g^2(x)] = f(G(x)) $\ or equivalently \ $(10')\ 2G(x) - xg(x) = g(x) F(G(x)) $ provides an isochronous  potential  verifying $G(0)=G'(0)=0, G''(0)=1$.\\

It is well known that the harmonic potential \ $ G=\frac{1}{2} x^2$\ is the only polynomial potential which is isochronous. Concerning the rational case, Chalykh and Veselov [2] proved the following 

\bigskip 

{\bf Proposition 3-8} \quad ([2]) {\it Under hypothesis ${(\cal H)}$ a rational potential $G(x)$ (which is not a polynomial) is isochronous if and only if} 
$$G(x) = \frac{1}{8 \alpha ^2} [\alpha x +1 - \frac{1}{\alpha x+ 1}]^2.$$

\bigskip

The proof we give below is different of that given by [2].

\bigskip

{\it Proof of Proposition 3-8} \quad Let us recall that for a involution $A$ we get\  $ G(A(x)) = G(x)$ \\ and \ $A(x) x < 0 $ for all $x \in [a,b]$. By Proposition 3-1 this potential is isochronous if
$$x - A(x) = 2\sqrt{2G(x)}\qquad for \quad all \qquad 0<x<b$$ or equivalently it verifies the functional equation $$G(x) = G(x-2\sqrt{2G(x)}\qquad for \quad all \qquad 0<x<b .$$
When the potential is rational it is required that $G(x)$ has to be the square of a rational function of the form 
$$G(x) = \frac {1}{2} (\frac {x P(x)}{Q(x)})^2$$ where $P(x)$ and $Q(x)$ are polynomials without common zeros.

This means the involution has to be meromorphic and has the same poles as $G(x)$. It can be written
$$A(x) = x - (\frac {2x P(x)}{Q(x)})= \frac{x Q(x) - 2x P(x)}{Q(x)}$$.

But, it is known (see for example Theorem (15.4) p. 296 of [6])  that the only meromorphic functions verifying $A^2(x) = x$  are $A(x) = L^{-1} (\frac{a}{L(x)})$ where $L(x)= a x + b$ a affine function since other meromorphic functions are not invertible. $A$ is then an homographic function\\ Thus, by hypothesis ${(\cal H)}$ and since $A(0)=0$ and $A'(0)=-1$ we get
$$A(x) = - \frac{x}{\alpha x +1}$$ implying $$\frac{2P(x)}{Q(x)}= 1+\frac{1}{\alpha x+1}=\frac{\alpha x+2}{\alpha x+1}.$$

\section{Parametrization of isochronous centers}
In this part we suppose that all functions are analytic. In particular, let us write
$$g(x)= x+\sum_{n\geq 2}a_n x^n\qquad and \qquad  G(x) = \frac{1}{2}x^2 + \sum_{n\geq 2} \frac {a_n}{n+1} x^{n+1},$$ and suppose that $r_0$ is the radius of convergence of these power series. By the Cauchy-Hadamard formula \ $\frac{1}{r_0}=lim_{n\rightarrow \infty}\mid a_n \mid^{1/n}.$\\
 
The purpose of this section is to highlight conditions so that an analytic potential $G(x)$ be isochronous. 
To that end let us write $$G(x) =  \frac{1}{2}x^2+ G_1(x)+G_2(x)$$ where the function \ $G_1(x)= \sum_{k\geq 2} \frac{a_{2k-2}}{2k-1}x^{2k-1}$ \ is odd and \\ $G_2(x) = \sum_{k\geq 2} \frac{a_{2k-1}}{2k}x^{2k}$\ is even.\\

We looking for general conditions on coefficients $a_n$ ensuring the isochronicity of the center $0$ of equation (1) \quad $\ddot x+g(x)=0.$\quad    Recall that only the harmonic potential $G(x) = \frac{1}{2}x^2$ is an even isochronous potential ($G_1(x)\equiv 0$). That means the odd coefficients $a_{2k+1}$ cannot be all zero when the isochronous potential is non-harmonic. In fact, we will prove little more.\\
   
 Starting from Theorem A we will show there are infinitely many necessary conditions verifying  by the coefficients in order the potential $G(x)$ be isochronous. More precisely, we state that the even coefficients $a_{2k}$ may be free and the odd coefficients $a_{2k+1}$ are polynomials with respect to $a_{2k}$.

\bigskip

{\bf Theorem C}\quad {\it Let the analytic potential $$G(x) = \frac{1}{2}x^2+\sum_{n\geq 3} \frac{a_{n-1}}{n}x^n=\frac{1}{2}x^2+G_1(x)+G_2(x)$$ of equation $$(1)\qquad \ddot x+g(x)=0.$$
When the equation (1) has an isochronous center at $0$ then the odd coefficients of the expansion of $g(x)$ can be expressed in terms of rational polynomials involving the even coefficients:  
$$ a_{2k+1}= f(a_{2k},a_{2k-2},...,a_2).  $$ In particular, when the potential $G(x)$ is isochronous then $G_2(x)\equiv 0$ is equivalent to $G_1(x)\equiv 0$, i.e. $G(x)= \frac{1}{2} x^2$ is harmonic.}\\

\bigskip

Thanks to {\it Maple} we are able to calculate the first terms :

$$a_3=\frac{10}{9}a_2^2,\qquad a_5=\frac{14}{5}a_2a_4-\frac{56}{27}a_2^4,$$ 

$$ a_7=\frac{ -592}{45}a_4a_2^3+\frac{848}{81}a_2^6+\frac{24}{7}a_2a_6+\frac{36}{25}a_4^2,$$

$$a_9=\frac{110}{27}a_2a_8-\frac{440}{21}a_2^3a_6+\frac{27808}{243}a_2^5a_4-\frac{536800}{6561}a_2^8-\frac{1144}{45}a_2^2a_4^2+\frac{22}{7}a_4a_6$$

\bigskip

${\it a_{11}}={\frac {52}{11}}\,{\it a_2}\,{\it a_{10}}+{\frac {57616}{
135}}\,{{\it a_2}}^{4}{{\it a_4}}^{2}-{\frac {2600}{81}}\,{{\it a_2
}}^{3}{\it a_8}+{\frac {125008}{567}}\,{{\it a_2}}^{5}{\it a_6}-{
\frac {4837664}{3645}}\,{{\it a_2}}^{7}{\it a_4}+{\frac {5631808}{
6561}}\,{{\it a_2}}^{10}-{\frac {2392}{125}}\,{\it a_2}\,{{\it a_4}
}^{3}-{\frac {7384}{105}}\,{{\it a_2}}^{2}{\it a_4}\,{\it a_6}+{
\frac {52}{15}}\,{\it a_4}\,{\it a_8}+{\frac {78}{49}}\,{{\it a_6}}
^{2}$\\

\bigskip

$a_{13} = -72\,{\it a_2}\,{\it a_6}\,{{\it a_4}}^{2}-{\frac {2632}{27}}\,{{
\it a_2}}^{2}{\it a_4}\,{\it a_8}+{\frac {38176}{27}}\,{{\it a_2}}
^{4}{\it a_4}\,{\it a_6}+{\frac {70}{13}}\,{\it a_2}\,{\it a_{12}}+{
\frac {42}{11}}\,{\it a_4}\,{\it a_{10}}+10/3\,{\it a_6}\,{\it a_8}-
{\frac {9430624}{1215}}\,{{\it a_4}}^{2}{{\it a_2}}^{6}+{\frac {
375769408}{19683}}\,{\it a_4}\,{{\it a_2}}^{9}+{\frac {166544}{225}}
\,{{\it a_4}}^{3}{{\it a_2}}^{3}-{\frac {920}{21}}\,{{\it a_2}}^{2}
{{\it a_6}}^{2}-{\frac {2190080}{729}}\,{{\it a_2}}^{7}{\it a_6}-{
\frac {14000}{297}}\,{{\it a_2}}^{3}{\it a_{10}}+{\frac {300944}{729}}
\,{{\it a_2}}^{5}{\it a_8}-{\frac {74681600}{6561}}\,{{\it a_2}}^{
12}-{\frac {616}{125}}\,{{\it a_4}}^{4}
$\\

\bigskip

 {\it Proof of Theorem C}\quad Let the analytic function 
 $$g(x)=x+\sum _{k=2}^{n}a_{{k}}{x}^{k}.$$ 
 By Theorem A the potential $G(x) = \frac{1}{2}x^2+\sum_{n\geq 3} \frac{a_{n-1}}{n}x^n$ is isochronous if and only if the following equality holds \ $(10)\quad  \frac{d}{dx}(\frac{G(x)}{g^2(x)})= f(G)$\quad where $f$ is an analytic function, set 
 
$$f(G) = b_0+b_1 G+b_2 G^2+ b_3 G^3+....$$ After replacing and equaling the two sides of (10) :
$$ \left( x+\sum _{k=2}^{n}a_{{k}}{x}^{k} \right) ^{-1}-(x^2+
\sum _{k=2}^{n}\frac{2 a_{{k}}}{k+1}{x}^{k+1}) \left( 1+\sum _{k=2}^{n}{\frac {a_{{
k}}{x}^{k}k}{x}} \right)  \left( x+\sum _{k=2}^{n}a_{{k}}{x}^{k}
 \right) ^{-3} =$$ $$\sum _{k=0}^{p}b_{{k}} \left( \frac{1}{2}\,{x}^{2}+\sum _{k=2}^{n}\frac{ a_{{k}}}{k+1}{x}^{k+1}) \right) ^{k}.$$

 Then we identify the analytic expansions of the two expressions. The unknown coefficients will then be determined by comparing powers in $x$:
$$b_0=-\frac{2}{3}a_2,\ a_3 =\frac{10}{9}a_2^2 \qquad b_1=11 a_2 a_3 -6 a_2^2-\frac{24}{5}a_4 $$ $$ -(10/3)\,a_{{5}}-(1/3)\,b_{{1}}a_{{2}}+{\frac {116}{15}}\,a_{{2}}a_{{4}}+
(14/3)\,{a_{{2}}}^{4}+4\,{a_{{3}}}^{2}-13\,a_{{3}}{a_{{2}}}^{2}=0,$$ 
$$-{\frac {30}{7}}\,a_{{6}}-\frac{1}{4}\,b_{{2}}-\frac{1}{4}\,b_{{1}}a_{{3}}+10\,a_{{2}}
a_{{5}}-{\frac {20}{3}}\,{a_{{2}}}^{5}+\frac{21}{2}\,a_{{4}}a_{{3}}-17\,a_{{4}
}{a_{{2}}}^{2}-{\frac {35}{2}}\,a_{{2}}{a_{{3}}}^{2}+25\,a_{{3}}{a_{{2
}}}^{3}=0,....$$

From these recursion formulae the coefficients can be easily determined. 
After eleminating $b_0, b_1, b_2,...$ we then deduce the expressions of coefficients \ $ a_3, a_5, a_7,...$.\\ 
The last part of Theorem C will proved by recurrence. Suppose that all coefficients $a_n=0$ for any $n < 2p$ where $p$ is a positive enteger and $b_0= b_1=b_2=....=b_{p-1}=0$. So, we may write 
$$g(x) = x+ a_{2p}x^{2p}+a_{2p+1}x^{2p+1}+a_{2p+2}x^{2p+2}+....$$ and 
$$G(x) = \frac{1}{2}x^2+\frac{a_{2p}}{2p+1}x^{2p+1}+\frac{a_{2p+1}}{2p+2}x^{2p+1}+\frac{a_{2p+2}}{2p+1}x^{2p+3}+...$$ Calculate
$$\frac{d}{dx}(\frac{G}{g^2})= -2p\,a_{{2\,p}}{x}^{2\,p-1}-a_{{2\,p+1}}{x}^{2\,p} \left( 2\,p+1
 \right) +a_{{2\,p}}{x}^{2\,p+1}+a_{{2\,p+1}}{x}^{2\,p+2} +.....$$ and equaling with $$ f(G)= b_p G^p+b_{p+1}G^{p+1}+...=\frac{b_p}{2^p}x^{2p}+ \frac{b_{p+1}}{2^{p+1}}x^{2p+2}+\frac{b_{p+2}}{2^{p+2}}x^{2p+4}+...$$
 Thus, necessarily $a_{{2\,p}}=0$.

\bigskip

{\bf Remarks 4-1}\\
Consider again Theorem B then equality 
$$(10)\qquad \frac{d}{dx}(G/g^2)= b_0 + b_1 G+ b_2 G^2+ b_3 G^3+ b_4 G^4+ ...$$ 
ensures the isochronicity of the potential $G$. \\ 
It is also possible to express coefficients $a_n$ of the analytic expansion of $g(x)$ (or its successive derivatives at $0$) in terms of $b_0, b_1, b_3,..$ \\
Indeed, 

$g'(0) = 1, g''(0) =-3 b_0, g'''(0) =15 b_0^2, g^{(4)}(0) = -105 b_0^3  - 30 b_1,\\ g^{(5)}(0) = 945 b_0^4  + 630 b_0 b_1, g^{(6)}(0) =-10395 b_0^5  - 11340 b_0^2  b_1 - 840 b_2,\\ g^{(7)}(0) =  207900 b_0^3 b_1 + 30240 b_0 b_2 + 135135 b_0^6  + 11340 b_1^3,...$\\

 $g$ may also be written 

$g(x) = x-\frac{3}{2}\,b_{{0}}{x}^{2}+\frac{5}{2}\,{b_{{0}}}^{2}{x}^{3}+ \left( -{\frac {5}{24
}}\,b_{{1}}-{\frac {35}{8}}\,{b_{{0}}}^{3} \right) {x}^{4}+{\frac {7}{
45}}\,b_{{0}} \left( {\frac {405}{8}}\,{b_{{0}}}^{3}+{\frac {45}{8}}\,
b_{{1}} \right) {x}^{5}+ \left( -{\frac {7}{120}}\,{\it b_2}-{\frac {
21}{8}}\,{b_{{0}}}^{2}b_{{1}}-{\frac {231}{16}}\,{b_{{0}}}^{5}
 \right) {x}^{6}+ \left( {\frac {55}{8}}\,{b_{{0}}}^{3}b_{{1}}+{\frac 
{429}{16}}\,{b_{{0}}}^{6}+\frac{3}{10}\,b_{{0}}b_{{2}}+\frac{1}{16}\,{b_{{1}}}^{2}
 \right) {x}^{7}+...$ \\
 
and the equation \ $\ddot x + g(x) = 0$\ has an isochronous center at $0$ for parameters values\ $b_0, b_1, b_2,...$\\

Thus, any isochronous potential $G$ is a multiparameters analytic function which may be expressed under the general form\\
 
$G(x) = \frac{1}{2}\,{x}^{2}-\frac{1}{2}\,b_{{0}}{x}^{3}+\frac{5}{8}\,{b_{{0}}}^{2}{x}^{4}+
 \left(-{\frac {1}{24}}\,b_{{1}}-{\frac {7}{8}}\,{b_{{0}}}^{3}
 \right) {x}^{5}+{\frac {7}{270}}\,b_{{0}} \left( {\frac {405}{8}}\,{b
_{{0}}}^{3}+{\frac {45}{8}}\,b_{{1}} \right) {x}^{6}+\frac{1}{7}\, \left( -{
\frac {7}{120}}\,{\it b_2}-{\frac {21}{8}}\,{b_{{0}}}^{2}b_{{1}}-{
\frac {231}{16}}\,{b_{{0}}}^{5} \right) {x}^{7}+\frac{1}{8}\, \left( {\frac {
55}{8}}\,{b_{{0}}}^{3}b_{{1}}+{\frac {429}{16}}\,{b_{{0}}}^{6}+\frac{3}{10}\,b
_{{0}}b_{{2}}+\frac{1}{16}\,{b_{{1}}}^{2} \right) {x}^{8}+...$\\

 A general expression of $a_{2p}$ as a polynomial of $b_0, b_1, b_2,...$ seems difficult to built. Nevertheless, it is possible to know the first coefficient $$a_{2p+1} = \frac{-1}{2^{p-1}}\frac{2p+1}{(2p)(2p-1)} b_{p-1}+...$$

\section {Application to the search of isochronous potentials}

In this part, we will see that  above results allow us to determinate families of isochronous potentials. More precisely, we then apply Theorems B to produce potentials with constant period. Some of them are new\\

{\large \bf A three-parameters family of isochronous potentials}\\ 

To be concrete, we derive at first a three-parameters family of potentials which appears little more general than the one given by Dorignac [6].  \\

Let us consider the following case

$$\frac {d}{dx}(\frac {G}{g^2})={\frac {\alpha}{ \left( 1+\beta G \right) ^{3/2}}} ,$$ where $\alpha $ and $\beta $ are real parameters such that \ $2\alpha^2 \leq \beta.$\\

Thanks to {\it Maple} a resolution of these equations yields
 
\begin{equation}G(x) =\frac {8\alpha^2+(\beta+2\alpha^2)(4\alpha x+\beta x^2)-(4\alpha^2+2\alpha\beta x)\sqrt{2(2+\beta x^2+4\alpha x)}}{2(\beta-2\alpha^2)^2}   \end{equation}
Then, the above potential is isochronous according to Theorem B. It may also be writen 
$$G(x) = \frac {1}{2}X^2$$ where 

$$ X=\frac{2\alpha+\beta x- \alpha \sqrt{2(2+\beta x^2+4\alpha x)}}{\beta-2\alpha^2}.$$
The function $h$ of Proposition 3-4 is
$$h(X) = \frac {\alpha \sqrt 2 X}{\sqrt{2+\beta X^2}}.$$
Recall that $X$ is such that $\frac{dX}{dx} > 0$ so that $X$ is a bijection of $x$. Moreover, $h(X)$ is odd function verifying $\mid h(X)\mid < 1 $ for any $X \in R$ since $2\alpha^2 \leq \beta.$ \\
 
 Applying scaling property of isochronous potentials, (Claim 2, Corollary of [6]). The potentials $G(x)$ and $\frac {1}{\gamma^2}G(\gamma x)$ have the same period. That means the following three-parameters potentials family is isochronous 
\begin{equation} G(x) = \frac {1}{2\gamma^2}X^2(\gamma x) = \frac{[2\alpha+\beta \gamma x- \alpha \sqrt{2(2+\beta \gamma^2 x^2+4\alpha \gamma x)}]^2}{2\gamma^2(\beta-2\alpha^2)^2}\end{equation}
  
The derivative of $G(x)$ is then

$$g(x)=\frac{2\alpha+\beta \gamma x- \alpha [\sqrt{2(2+\beta \gamma^2 x^2+4\alpha \gamma x)}](\beta \gamma -\frac{\alpha(4\beta\gamma^2 x+8\alpha\gamma )}{ \sqrt{2(2+\beta \gamma^2 x^2+4\alpha \gamma x)}}}{\gamma^2(\beta-2\alpha^2)^2} 
 $$
 The involution $A$ defined by $G(A(x)=G(x), \ A(x)x<0$\  may be written
 $$A(x) = x - 2  \frac{[2\alpha+\beta \gamma x- \alpha \sqrt{2(2+\beta \gamma^2 x^2+4\alpha \gamma x)}]}{\gamma(\beta-2\alpha^2)}$$

As special cases we may derive the following \\

{\bf 1 -} For any $ \gamma \neq 0$ and $\beta = 2\alpha$ we obtains the two-parameters family of isochronous potentials (see (24) of [6]) introduced  for the first time by Stillinger and Stillinger [12]. Indeed, after replacing in (21) $\beta $ by $2 \alpha$ and simplifying by $2 \alpha$ it yields 

\begin{equation} G(x)  = \frac{[1+ \gamma x-  \sqrt{1+\alpha \gamma^2 x^2+2\alpha \gamma x)}]^2}{2(1-\alpha)^2}\end{equation}

{\bf 2 -} The case $\alpha=\beta=0$ and $\gamma = 1$ yields the harmonic potential : \ $G(x) = \frac {1}{2}x^2.$\\

{\bf 3 -} The case $\beta=0$ and $\gamma = 1$ gives the Urabe potential (see Corollary 2-2) :  $$G(x) = \frac {4}{\alpha^2} - \frac {2}{\alpha}(x + 2\frac {\sqrt{1-\alpha x}}{\alpha}).$$

{\bf 4 -} The case $2\alpha=\beta$ and $\gamma = 1$ yields the Bolotin-Mc Kay potential (see [6]).
Indeed, $$\frac {d}{dx}(\frac {G}{g^2}) = {\frac {\alpha}{ \left( 1+2 \alpha G \right) ^{3/2}}}$$ implies
$$G(x)  = \frac{[1+ x-  \sqrt{1+\alpha x^2+2\alpha x)}]^2}{2(1-\alpha)^2}$$
We then deduce the function $h$ of Proposition 3-4 
$$h(X)= \,{\frac {\alpha X}{\sqrt {1+\,{X}^{2}{\alpha}}}}$$

{\bf 5 -} The case $2\alpha^2=\beta$ and $\gamma = 1$ yields the isotonic potential (see [6])
Thanks to {\it Maple} a resolution of 
$$\frac {d}{dx}(\frac {G}{g^2})= {\frac {\alpha}{ \left( 1+ \alpha^2 G \right) ^{3/2}}}$$
gives
$$G={\frac {1}{4}}\,{\frac {x^2(2+\alpha x)^2}{(1+\alpha x)^2}}.$$

The function $h$ of Proposition 3-4 is
$$h(X) = {\frac {\alpha X}{\sqrt {1+\,{X^2}{\alpha}^{2}}}}$$  and its integral is

$$H(X)={\frac {\sqrt {1+\,{X}^{2}{\alpha}^{2}}}{\alpha}} - \frac{1}{\alpha} = x - X$$

Thus $${{\frac {1+\,{X}^{2}{\alpha}^{2}}{{\alpha}^{2}}}= \left( x-X+
{{\alpha}^{-1}} \right) ^{2}}$$ 

We then obtain by another way
$$G={\frac {1}{4}}\,{\frac {x^2(2+\alpha x)^2}{(1+\alpha x)^2}}= \frac{1}{8\alpha^2}[\alpha x+1 -\frac{1}{\alpha x+1}]2$$

\bigskip

{\large \bf Others isochronous potentials}\\

We will give others two-parameters families of potentials with constant period which seems to be new in the literature.\\ 

{\bf 1 -}  
Let us consider  
$$ \frac {d}{dx}(\frac {G}{g^2})= {\frac {\alpha}{ \left( 1+2 \beta^2 G \right) ^{5/2}}}.$$

Then the function $h$ of Propositon 3-4 may be calculated
$$h(X)=\alpha \left( 1/3\,{\frac {X}{ \left( 1+{X}^{2}{\beta}^{2} \right) ^{3/2}}}+2/3
\,{\frac {X}{\sqrt {1+{X}^{2}{\beta}^{2}}}} \right)$$
and its integral is
$$H(X)=1/3\,{\frac {\alpha \left( 1+2\,{X}^{2}{\beta}^{2} \right) }{{\beta}^{2}\sqrt {1+{X
}^{2}{\beta}^{2}}}} - \frac {\alpha}{3\beta ^2}$$

Taking for example\ $\beta^2=2\alpha/3$\ and thanks to {\it Maple} we find another isochronous potential. \\
 Indeed, since $x =X+H(X) $ and solving $H(X)^2= (x - X)^2$ \ we get\\

$
X(x)  =  (1/24)\,{\frac { [u(x)]^{1/3}}{\alpha \left( 3+4\,x \alpha \right) }}+(1/24)\,{\frac { \left( 16
\,{x}^{2}{\alpha}^{2}+24\,x \alpha-423 \right)  \left( 3+4\,x \alpha \right) }{\alpha [u(x)]^{1/3}}+(1/24)
\,{\frac {3+4\,x \alpha}{\alpha}}}
$\\

where\\
  $ u(x)= -5751+31536\,x\alpha+21600\,{x}^{2}{\alpha}^{2
}+768\,{x}^{3}{\alpha}^{3}+256\,{x}^{4}{\alpha}^{4}$ \\ $+72\,\sqrt {3}\sqrt {1024\,{x
}^{6}{\alpha}^{6}+4608\,{x}^{5}{\alpha}^{5}+26496\,{x}^{4}{\alpha}^{4}+62208\,{x}^{3}
{\alpha}^{3}+101412\,{x}^{2}{\alpha}^{2}+86022\,x\alpha+45927} $\\

So, according to Theorem A and in using scaling properties the potential \\

$G(x) = (\frac{1}{2})(\frac{1}{24})^2 [{\frac { [u(\gamma x)]^{1/3}}{\alpha \left( 3+4\,\gamma x \alpha \right) }}+{\frac { \left( 16
\,{\gamma^2 x}^{2}{\alpha}^{2}+24\,\gamma x \alpha-423 \right)  \left( 3+4\,\gamma x \alpha \right) }{\alpha [u(\gamma x)]^{1/3}}+
{\frac {3+4\,\gamma x \alpha}{\alpha}}}]^2$\\

is isochronous.\\

{\bf 2 -}  Let $$ \frac {d}{dx}(\frac {G}{g^2})= {\frac {2\alpha + 2 \alpha^3 G}{ \left( 1+2 \alpha^2 G \right) ^{5/2}}}.$$
The function $h$ is then
$$h(X) = \frac{2\alpha X + \alpha^3 X^3}{(1+\alpha^2 X^2)}$$
and its integral has the following simple form
$$H(X) = \frac{\alpha X^2}{\sqrt{1+\alpha^2 X^2}}.$$

Thanks to {\it Maple} we find 
$$X = \frac{1}{6 \alpha^2 x}v(x) + \frac{(-10\alpha^2 x^2 + 1+\alpha^4 x^4)}{6\alpha^2 x v(x)}+ \frac{1+\alpha^2 x^2}{6\alpha^2 x}$$
where $$v(x) = [-15\alpha^2 x^2+39\alpha^4 x^4+1+\alpha^6 x^6 +\alpha^3 x^3\sqrt{(-3+33\alpha^2 x^2+3\alpha^4 x^4)}]^{1/3}$$

\bigskip
{\bf Aknowledgments}\\ {\it I would like to thank Jean-Marie Strelcyn. The final form of this paper owes him very much.}\\

\bigskip

{\bf References}

\bigskip

[1]\ F. Calogero  \quad {\it Isochronous systems}\quad  Oxford University Press, Oxford, (2008).  

\smallskip 

[2] \ O. Chalykh and A. Veselov \quad {\it  A remark on rational isochronous potentials} \quad J. Nonlinear Math. Phys. {\bf 12-1}, p. 179-183, (2005). 

\smallskip        

[3]\  R. Chouikha and F. Cuvelier \quad {\it Remarks on some monotonicity conditions for the period function}\quad  Applic. Math., {\bf 26}, no. 3, p. 243-252, (1999).

\smallskip 
[4]\  R. Chouikha \quad {\it Monotonicity of the period function for some planar differential systems, I. Conservative and quadratic systems} \quad     Applic. Math., {\bf 32}, no. 3, p. 305-325, (2005). 
                        
\smallskip 
[5]  \ S.N. Chow and D. Wang \quad {\it On the monotonicity of the period function of some second order equations}\quad Casopis Pest. Mat. {\bf 111}, p. 14-25, (1986).

\smallskip
[6] \ J. Dorignac \quad {\it On the quantum spectrum of isochronous potentials}\quad J. Phys; A: Math. Gen., {\bf 38}, p. 6183-6210 (2005).

\smallskip
[7] \ I. Koukles and N. Piskounov \quad {\it Sur les vibrations tautochrones dans les systèmes conservatifs et non conservatifs}. C. R. Acad. Sci., URSS, vol XVII, n°9, p. 417-475, (1937).

\smallskip
[8]\ M. Kuczma \quad {\it Functional equations in a single variable} \quad Monografie Matematyczne, Tom 46, Warsaw, (1968).

\smallskip
[9]\ L.D. Landau E.M. Lifschitz {\it Mechanics, Course of Theorical Physics} Vol 1, Pergamon Press, Oxford, (1960).                     

\smallskip
[10] \ F. Rothe \quad {\it Remarks on periods of planar Hamiltonian systems.} \quad SIAM J. Math. Anal., {\bf 24}, p.129-154, (1993).

\smallskip
[11] \ R. Schaaf \quad {\it A class of Hamiltonian systems with increasing periods} \quad J. Reine Angew. Math., {\bf 363}, p. 96-109, (1985). 

\smallskip
[12]\ FH. Stillinger and DK. Stillinger \quad {\it Pseudoharmonic oscillators and inadequacy of semiclassical quantization}\quad J. Phys. Chim., {\bf 93}, 6890, (1989).

\smallskip
[13]\  M. Urabe \quad  {\it The potential force yielding a periodic motion whose period is an arbitrary continuous function of the amplitude of the velocity}\ Arch. Ration. Mech. Anal.,{\bf 11}, p.27-33, (1962). 

\end{document}